\documentclass{elsart}
\usepackage{amsfonts,amssymb,amsmath}
\begin{document}
\begin{frontmatter}
\title{Involutions on graded matrix algebras}    
\author{Yuri Bahturin}\address{Department of Mathematics and Statistics\\Memorial
University of Newfoundland\\ St. John's, NL, A1C5S7, Canada\\
 and \\ Department of
Algebra, Faculty of Mathematics and Mechanics\\Moscow State University\\Moscow, 119992, Russia}\ead{yuri@math.mun.ca}\thanks{Work is partially supported by NSERC grant \#
227060-04 and URP grant, Memorial University of Newfoundland}
\author{Mikhail Zaicev}\address{Department of Algebra,
Faculty of Mathematics and Mechanics\\Moscow State University\\Moscow, 119992, 
Russia}\ead{zaicev@mech.math.msu.su}\thanks{Work is partially supported by RFBR, grant 06-01-00485a, and SSH-5666.2006.1}                   

\end{frontmatter}

\textwidth 140mm 
\textheight 210mm
\topmargin -0mm
\oddsidemargin 10mm
\evensidemargin 5mm
\newcommand{\pp}{\noindent {\em Proof. }}
\newtheorem{theorem}{Theorem}
\newtheorem{lemma}{Lemma}
\newtheorem{proposition}{Proposition}
\newtheorem{remark}{Remark}
\newcommand{\bee}[1]{\begin{equation}\label{#1}}
\newcommand{\beq}[1]{\begin{eqnarray}\label{#1}}
\newcommand{\ben}{\begin{eqnarray*}}
\newcommand{\ene}{\end{equation}}
\newcommand{\eqe}{\end{eqnarray}}
\newcommand{\eqn}{\end{eqnarray*}}
\newcommand{\ld}{\ldots}
\newcommand{\diag}{\mathrm{diag}}
\newcommand{\cd}{\cdots}
\newcommand{\tr}[1]{\: ^t\!#1}
\newcommand{\rsl}{\mathrm{sl}}
\newcommand{\rgl}{\mathrm{gl}}
\newcommand{\wh}[1]{\widehat{#1}}
\newcommand{\wg}{\widehat{G}}
\newcommand{\wR}{\widetilde{R}}
\newcommand{\vp}{\varphi}
\newcommand{\bx}{\hfill$\Box$}
\newcommand{\cj}[2]{#1^{-1}\!#2#1}
\newcommand{\iv}[1]{#1^{-1}\!}
\newcommand{\ve}{\varepsilon}
\newcommand{\xba}{X_{\bar{a}}}
\newcommand{\xbb}{X_{\bar{b}}}
\newcommand{\xbab}{X_{\bar{a}\bar{b}}}
\newcommand{\xbe}{X_{\bar{e}}}
\newcommand{\xbg}{X_{\bar{g}}}
\newcommand{\ba}{\bar{a}}
\newcommand{\bz}{\bar{0}}
\newcommand{\bw}{\bar{1}}
\newcommand{\be}{\bar{e}}
\newcommand{\bg}{\bar{g}}
\newcommand{\bG}{\bar{G}}
\newcommand{\bh}{\bar{h}}
\newcommand{\bp}{\bar{p}}
\newcommand{\bq}{\bar{q}}
\newcommand{\bP}{\bar{P}}
\newcommand{\bQ}{\bar{Q}}
\newcommand{\bi}{\bar{i}}
\newcommand{\bj}{\bar{j}}
\newcommand{\bu}{\bar{u}}
\newcommand{\by}{\bar{y}}
\newcommand{\bv}{\bar{v}}
\newcommand{\epf}{\hfill\qed}
\newcommand{\su}[1]{\mathrm{Supp}\,#1}
\newcommand{\Sp}[1]{\mathrm{Span}\,\{#1\}}
\newcommand{\twm}[4]{\left(\begin{array}{cc}#1&#2\\#3&#4\end{array}\right)}


\begin{abstract}
In this paper we describe graded automorphisms and antiautomorphisms of finite order on matrix algebras endowed with a group gradings by a finite abelian group over an arbitrary algebraically closed field of charcteristic different from 2.
\hfill\break 
 \hfill\break 
{\sl Mathematics Subject Index 2000:} 17B20, 17B40, 16W10, 16W50 
 \end{abstract}



\section{Introduction}\label{s0}

 This paper is devoted to the correction of an error in the paper \cite{BShZ} in which the classification of involution gradings on matrix algebras was derived from the fact that in the decomposition of a graded matrix algebra as the tensor product of an elementary and a fine component, these components remain invariant under the involution.

\section{Some notation and simple facts}\label{s1}

 Let $F$ be an arbitrary field, $A$ a not necessarily associative algebra over an $F$ and $G$ a group. We say
that  $A$ is a $G$-graded algebra, if there is a vector space sum
decomposition
\bee{e0001}
A=\bigoplus_{g\in G} A_g,
\ene
such that
\bee{e0002}
A_gA_h\subset A_{gh}\mbox{ for all }g,h\in G.
\ene

A subspace $V\subset A$ is called {\em graded} (or {\it homogeneous})
if $V=\oplus_{g\in G} (V\cap A_g)$. An element $a\in R$ is called
{\it homogeneous of degree} $g$ if $a\in A_g$. We also write $\deg a=g$. The {\it support} of the
$G$-grading is a
subset
$$
\mathrm{Supp}\: A=\{g\in G|A_g\ne 0\}.
$$

Suppose now that $F$ is of characteristic different from $2$. If $A$ is an associative algebra with involution $\ast$ and, in addition to (\ref{e0002}),  one has
\bee{e0004}
(A_g)^\ast= A_g\mbox{ for all }g\in G.
\ene
then we say that (\ref{e0001}) is an \textit{involution preserving} grading or simply an \textit{involution} grading. In this case, given a graded subspace $B\subset A$ we set
\bee{esym}
H(B,\ast)=\{ b\in B\,|\,b^*=b\},\mbox{ the set of symmetric elements of }B
\ene
and
\bee{essym}
K(B,\ast)=\{ b\in B\,|\,b^*=-b\},\mbox{ the set of skew-symmetric elements of }B.
\ene
If $B$ is an associative subalgebra of $A$ then $B^{(-)}$ is a Lie subalgebra of $A$, that is, with respect to $[x,y]=xy-yx$ while $B^{(+)}$ is a Jordan subalgebra of $A$, that is, with respect to $x\circ y=xy+yx$. We always have $B=B^{(-)}\oplus B^{(+)}$.

\section{Reminder: Group gradings on matrix algebras}\label{stwotypes}

Below we briefly recall the results of \cite{BSZ}, where the full description of abelian group gradings on the full matrix algebra has been given.

\bigskip

A grading $R=\oplus_{g\in G}R_g$ on
the matrix algebra $R=M_n(F)$ is called {\it elementary} if there exists an $n$-tuple
$(g_1,\ldots,g_n)\in G^n$ such that the matrix units $E_{ij}, 1\le
i,j\le n$ are homogeneous and $ E_{ij}\in R_g\iff g=g_i^{-1}g_j. $

A grading is called {\it fine} if $\dim R_g=1$ for any $g\in
\mathrm{Supp}\: R$. A particular case of fine gradings is the so-called
$\varepsilon$-grading where $\varepsilon$ is $n^\mathrm{th}$ primitive root
of $1$. Let $G=\langle a\rangle_n\times\langle b\rangle_n$ be the direct product of two cyclic groups of order
$n$ and

\bee{e00}
X_a = \left(\begin{array}{cccc} \varepsilon^{n-1} & 0 & ... & 0
\\ 0 & \varepsilon^{n-2} & ... & 0
\\  \cdots& \cdots & \cdots& \cdots
\\ 0 & 0 & ... & 1
\end{array}\right)~,~~
X_b = \left(\begin{array}{cccc}
  0 & 1 &  ... & 0
\\ \cdots & \cdots   &    \cdots &  \cdots
\\ 0 & 0 &  ... & 1
\\ 1 & 0 &  ... & 0
\end{array}\right)~.
\ene
Then
\begin{equation}\label{a1}
X_aX_bX_a^{-1}=\varepsilon X_b~,~~X_a^n=X_b^n=I
\end{equation}
and all $X_a^iX_b^j, 1\le i,j\le n$, are linearly independent.
Clearly, the elements $X_a^iX_b^j,\, i,j=1,\ldots, n$, form a basis
of $R$ and all the products of these basis elements are uniquely
defined by (\ref{a1}).

       Now for any $g\in G, g=a^ib^j$, we set $X_g=X_a^iX_b^j$ and denote by $R_g$
a one-dimensional subspace
\bee{aa1}
R_g=\langle X_a^iX_b^j\rangle.
\ene
Then from (\ref{a1}) it follows that $R=\oplus_{g\in G}R_g$ is a
$G$-grading on $M_n(F)$ which is called an $\varepsilon$-grading.

Now let $R=M_n(F)$ be the full matrix algebra over $F$ graded by an abelian
group $G$. The folowing result has been proved in \cite[Section 4, Theorems 5, 6]{BSZ} and \cite[Subsection 2.2, Theorem 6, Subsection 2.3, Theorem 8]{surgrad}.

\begin{theorem}\label{tbsz} Let $F$ be an algebraically closed field of characteristic zero. Then as a $G$-graded algebra $R$ is isomorphic to the
tensor product
$$
R^{(0)}\otimes R^{(1)}\otimes\cd\otimes R^{(k)}
$$
where $R^{(0)}=M_{n_0}(F)$ has an elementary $G$-grading, $\mathrm{Supp}\:
R^{(0)}=S$ is a
finite subset of $\,$ $G$, $R^{(i)}=M_{n_i}(F)$ has the $\varepsilon_i$ grading,
$\varepsilon_i$ being a primitive $n_i^\mathrm{th}$ root of $1$,
$\mathrm{Supp}\; R^{(i)}=H_i\cong \mathbb{Z}_{n_i}\times
\mathbb{Z}_{n_i}, i=1,\ld,k$. Also $H=H_1\cd H_k\cong H_1\times\cd\times
H_k$ and $S\cap H=\{e\}$ in $G$.
\end{theorem}

\begin{remark}\label{r101}
It follows from a very general lemma in \cite{BSZ} that  the support $T$ of a fine grading on $R=M_n$  is a subgroup of the grading group $G$. Thus we have $R=\bigoplus_{t\in T}R_t$ and $R_t=\langle X_t\rangle$, for a nondegenerate matrix $X_t$. Let us also recall that the product in $R=M_n(F)$ with fine grading as above is defined by a \textit{bicharacter} $\alpha:T\times T\rightarrow F^\ast$ as follows: $X_tX_u=\alpha(t,u)X_{tu}$, for any $t,u\in T$. The commutation relations in $R$ take the form $X_tX_u=\beta(t,u)X_uX_t$ where $\beta(t,u)=\alpha(t,u)/\alpha(u,t)$ is a
\textit{skew-symmetric bicharacter} on $T$ (see \cite{BMZ}).
\end{remark}

Let us recall that any involution $\ast$ of $R=M_n$ can always be written as

\bee{einv}
X^\ast= \Phi^{-1}( \tr{X})\Phi
\ene
where $\Phi$ is a nondegenerate matrix which is either symmetric or skew-sym\-metric and $X\mapsto\, \tr{X}$ is the ordinary transpose map. In the case where $\Phi$ is symmetric we call $\ast$ a transpose involution. If $\Phi$ is skew-symmetric $\ast$ is called a symplectic involution. Before we formulate the theorem describing involution gradings on $M_n$ in the case where the elementary and fine components are invariant under the involution, we need three (slightly modified) lemmas from \cite{BShZ}. The general restriction in \cite{BShZ} zero characteristic was not used in the proof of these particular lemmas. The first two deal with elementary involution gradings while the last with certain fine involution gradings. If $R$ has an involution $\ast$ then by $R^{(\pm)}$ we denote the space of symmetric (respectively skew-symmetric) matrices in $R$ under $\ast$.

The next lemma handles the case of an elementary grading compatible with an involution defined by a symmetric non-degenerate bilinear form.

\begin{lemma}\label{L8'}
Let $R=M_n(F)$, $n$ a natural number, be a matrix algebra with involution $\ast$
defined by a symmetric non-degenerate bilinear form. Let $G$ be an
abelian group and let $R$ be equipped with an elementary
involution $G$-grading defined by an $n$-tuple $(g_1,\ld,g_n)$. Then
$g_1^2=\ld=g_m^2=g_{m+1}g_{m+l+1}=\ld=g_{m+l}g_{m+2l}$ for some
$0\le l\le\frac{n}{2}$ and $m+2l=n$. The involution $\ast$  acts as
$X^\ast =(\Phi^{-1})\:^tX\Phi$ where
$$
\Phi=\left(\begin{array}{ccc}   I_m   &   0     &  0   \cr
                  0    &   0     &  I_l \cr
                  0    &   I_l   &  0   \cr\end{array}\right),
$$
where $I_s$ is the $s\times s$ identity matrix. Moreover, $R^{(-)}$
consists of all matrices of the type
\bee{BD1}
      \left(\begin{array}{ccc}  P    &   S     &   T     \cr
               -^t T   &   A     &   B     \cr
               -^t S   &   C     & -^t A    \cr\end{array}\right),
\ene
where $^t P=-P,~ ^t B=-B,~ ^t C=-C$ and
$$
P\in M_m(F),~ A,B,C,D\in M_l(F),~ S,T\in M_{m\times l}(F)
$$
while $R^{(+)}$ consists of all matrices of the type
\bee{BD2}
      \left(\begin{array}{ccc}  P    &   S     &   T    \cr
                ^t T   &   A     &   B    \cr
                ^t S   &   C     &  ^t A   \cr\end{array}\right),
\ene
where $^t P= P,~ ^t B= B,~ ^t C= C$ and
$$
P\in M_m(F),~ A,B,C,D\in M_l(F),~ S,T\in M_{m\times l}(F).
$$
\end{lemma}

The next lemma deals with the case of an elementary grading compatible with an involution defined by a skew-symmetric non-degenerate bilinear form.

\begin{lemma}\label{L8}
Let $R=M_n(F)$, $n=2k$, be the matrix algebra with involution $\ast$
defined by a skew - symmetric non - degenerate bilinear form. Let $G$ be
an abelian group and let $R$ be equipped with an elementary
involution $G$-grading defined by an $n$-tuple $(g_1,\ld,g_n)$. Then $g_{1}g_{k+1}=\ld=g_{k}g_{2k}$,
the involution $\ast$ acts as $X^\ast =(\Phi^{-1})\:^tX\Phi$ where
$$
\Phi=\left(\begin{array}{cc}  0   &   I     \cr
                -I    & 0    \cr\end{array}\right),
$$
$I$ is  the $k\times k$ identity matrix, $R^{(-)}$ consists of all
matrices of the type
\bee{C1}
      \left(\begin{array}{cc}  A    &   B     \cr
                 C    & -^t A    \cr\end{array}\right),~
A,B,C,\in M_k(F),~ ^t B=B,^t C=C
\ene
while $R^{(+)}$ consists of all matrices of the type
\bee{C2}
      \left(\begin{array}{cc}  A    &   B     \cr
                 C    &  ^t A    \cr\end{array}\right),~
A,B,C,\in M_k(F),~ ^t B=-B,^t C=-C.
\ene
\end{lemma}
\begin{lemma}\label{L6}
Let $R=M_2(F)$ be a $2\times 2$ matrix algebra endowed with an involution $*:R\to R$ corresponding to a symmetric or skew-symmetric
non-degenerate bilinear form with the matrix $\Phi$. The $(-1)$-grading
of $M_2$ by $G=\langle a\rangle_2\times \langle b\rangle_2$ is an involution grading if and only if one of the following
holds:
\begin{itemize}
\item[(1)] $\Phi$ is skew-symmetric,
$$
\Phi=\twm{0}{1}{-1}{0},
\quad
K(R,\ast)=\Sp{X_a,X_b,X_{ab}},
\quad
H(R,\ast)=\Sp{X_e};
$$

\item[(2)] $\Phi$ is symmetric,
$$
\Phi=\twm{0}{1}{1}{0},
\quad
K(R,\ast)=\Sp{X_a},
\quad
H(R,\ast)=\Sp{X_e,X_b,X_{ab}};
$$
\item[(3)] $\Phi$ is symmetric,
$$
\Phi=\twm{1}{0}{0}{1},
\quad
K(R,\ast)=\Sp{X_{ab}},
\quad
H(R,\ast)=\Sp{X_e,X_a,X_b};
$$

\item[(4)] $\Phi$ is symmetric,
$$
\Phi=\twm{1}{0}{0}{-1},
\quad
K(R,\ast)=\Sp{X_b},
\quad
H(R,\ast)=\Sp{X_e,X_a,X_{ab}};
$$\end{itemize}
\end{lemma}

Notice that the involution in each case is already defined, say, in Case $\mathrm(1)$ one has
$$(\alpha X_e+\beta X_a+\gamma X_b+\delta X_{ab})^\ast=\alpha X_e-\beta X_a-\gamma X_b-\delta X_{ab}.$$

\begin{remark}\label{r102}If $R=M_n$ has a fine grading by a group $G$ with support an elementary abelian $2$-subgroup $T$ then it is immediate from the previous lemma and a Remark \ref{r101} after Theorem \ref{tbsz} that $R$ has a basis $\{ X_t\,\vert\, t\in T\}$ such that $X_tX_u=\alpha(t,u)X_{tu}$ where $\alpha(t,u)=\pm 1$ and for each $u\in T$ we have $\iv{X_u}=\tr{X_u}=\alpha(u,u)X_u$.
\end{remark}

\bigskip

We can now formulate the most general result available earlier, which describes gradings on a matrix algebra with involution (a weaker form of \cite[Theorem 2]{BShZ}, which is not true). With our additional assumption that the involution respects the fine and the elementary components of the grading, the proof of \cite{BShZ} works without changes. We remark here that this condition is always satisfied provided that the Sylow 2-subgroup of $G$ is cyclic.

\begin{theorem}\label{inv}
Let $R=M_n(F)=\oplus_{g\in G}R_g$ be a matrix algebra over an
algebraically closed field of characteristic zero graded by the
group $G$ and $\mathrm{Supp}\: R$ generates $G$. Suppose that
$*: R\rightarrow R$ is a graded involution. Then $G$ is abelian,
and $R$ as a $G$-graded algebra  is
isomorphic to the tensor product $R^{(0)}\otimes
R^{(1)}\otimes\cdots\otimes R^{(k)}$ of a matrix subalgebra $R^{(0)}$ with elementary grading and $R^{(1)}\otimes\cdots\otimes R^{(k)}$ a matrix subalgebra with fine grading. Suppose further that both these subalgebras are invariant under the involution. Then $n=2^km$ and
\begin{itemize}
\item[(1)]
$R^{(0)}=M_m(F)$ is as in Lemma \ref{L8} if $\ast$ is symplectic on
$R^{(0)}$ or as in Lemma \ref{L8'} if $\ast$ is transpose on $R_0$;
\item[(2)]
$R^{(1)}\otimes\cdots\otimes R^{(k)}$ is a $T=T_1\times\ldots\times
T_k$-graded algebra and any $R^{(i)},1\le i\le k$, is $T_i\cong
{\mathbb Z}_2\times{\mathbb Z}_2$-graded algebra as in Lemma
\ref{L6}.
\item[(3)] A graded basis of $R$ is formed by the elements $Y\otimes X_{t_1}\otimes\cdots\otimes X_{t_k}$, where $Y$ is an element of a graded basis of $R^{(0)}$ and the elements $X_{t_i}$ are of the type (\ref{aa1}), with $n=2$, $t_i\in T_i$. The involution on these elements is given canonically by
$$(Y\otimes X_{t_1}\otimes\cdots\otimes X_{t_k})^\ast=Y^\ast\otimes X_{t_1}^\ast\otimes\cdots\otimes X_{t_k}^\ast=\mathrm{sgn}(t)(Y^\ast\otimes X_{t_1}\otimes\cdots\otimes X_{t_k}),$$
where $Y\in R^{(0)}$, $X_{t_i}$ are the elements of the basis of the canonical $(-1)$-grading of $M_2$, $i=1,\ldots,k$, $t=t_1\cdots t_k\in T$, $\mathrm{sgn}(t)=\pm 1$, depending on the cases in Lemma \ref{L6}.

\end{itemize}
\end{theorem}

In the next two sections we describe the antiautomorphisms of graded matrix algebras in the general case, including that now we only assume that the base field is algebraically closed of characteristic different from 2.

\section{Antiautomorphisms of graded matrix algebras}\label{sAA}

We start this section with a result about the structure of fine gradings of $R=M_n$ compatible with an antiautomorphism. This result is a far going generalization of \cite[Lemma 2]{BShZ}.
Any antiautomorphism $\vp$ of $R=M_n$ can always be written as

\bee{einv1}
\vp\ast X= \Phi^{-1}( \tr{X})\Phi
\ene
where $\Phi$ is a nondegenerate matrix and $X\mapsto\, \tr{X}$ is the ordinary transpose map. It is well-known that $\vp$ is an involution if and only if $\Phi$ is either symmetric or skew-symmetric. Recall that in the case where $\Phi$ is symmetric $\vp$ is called a \textit{transpose involution} and if $\Phi$ is skew-symmetric then $\vp$ is called a \textit{symplectic involution}.

\begin{lemma}\label{lfin1}
Let $R=M_n(F)=\bigoplus_{t\in T}R_t$ be the $n\times n$-matrix algebra with an
$\ve$-grading, $T=\langle a\rangle_n\times \langle b\rangle_n$. Let also $\vp:
R\rightarrow R$ be an antiautomorphism of $R$ defined by $\vp\ast X=\cj{\Phi}{\,\tr{X}}$. If
$\vp*R_t=R_t$ for all $t\in T$ then $n=2$, $\Phi$ coincides with the scalar multiple one of the matrices $I$,
$X_a$, $X_b$ or $X_{ab}$ (see (\ref{e00})).
\end{lemma}

\pp First we consider the $\vp$-action on $X_a$. Since $R_a$ is stable under $\vp$,
$$
\Phi^{-1~t}X_a\Phi=\Phi^{-1}X_a\Phi=\alpha X_a
$$
for some scalar $\alpha\ne 0$. Then \bee{f1} X_a\Phi X_a^{-1}=\alpha\Phi. \ene Since $X_a^n=I$, we obtain $\alpha^n=1$, so that $\alpha=\ve^j$ for
some $0\le j\le n-1$.

Denote by $P$ the linear span of $I,X_a,\ld,X_a^{n-1}$. Then $R=P\oplus X_bP\oplus \cd
\oplus X_b^{n-1}P$ as a vector space and the conjugation by $X_a$ acts on $X_b^iP$ as the
multiplication by $\ve^i$. In particular, all eigenvectors with eigenvalue $\ve^j$ are in
$X_b^jP$. It follows that $\Phi\in X_b^jP$, that is, $\Phi= X_b^jQ$ for some $Q\in P$.

Now we consider the action of $\vp$ on $X_b$:
$$
\vp* X_b=\Phi^{-1~t}X_b\Phi=\Phi^{-1}X_b^{-1}\Phi=\gamma X_b,
$$
that is, $X_b\Phi X_b = \mu\Phi$ with $\mu=\gamma^{-1}\ne 0$. If we write
$Q=\sum\alpha_i X_a^i$ then \bee{f2} X_b\Phi X_b = X_b^j\sum_{i}\alpha_iX_b X_a^i X_b =
X_b^j\sum_{i}\alpha_i' X_a^i X_b^2=\mu\Phi=\mu X_b^j\sum_{i}\alpha_iX_a^i, \ene In this
case $X_b^j\sum_{i}\alpha_i' X_a^i X_b^2=\mu X_b^j\sum_{i}\alpha_iX_a^i$ where the
scalars $\alpha_i'$ can be explicitly computed using (\ref{f1}). Since the degrees in
$X_a,X_b$ define the degrees in the $T$-grading, we can see that (\ref{f2}) immediately
implies $X_b^2=I$, i.e. $n=2$.

As we have shown before, (\ref{f1}) implies $\Phi=X_b^j Q$ with $Q=\alpha_0 I+\alpha_1
X_a$. Since $n=2$, the argument following (\ref{f1}) applies if we change $a$ and $b$
places so that $\Phi=X_a^k(\beta_0 I + \beta_1 X_b)$. Comparing these two expressions we
obtain that $\Phi$ must be one of $I$, $X_a$, $X_b$, or $X_{ab}$, up to a scalar
multiple.  \bx

Now we make few remarks about the structure of elementary gradings on $M_n(F)$. Recall that a grading
$M_n=R=\oplus_{g\in G}R_g$ is elementary if there exists an $n$-tuple $\tau =
(g_1,\ldots,g_n)\in G^n$ such that the matrix units $E_{ij}, 1\le i,j\le n$ are
homogeneous and $ E_{ij}\in R_g\iff g=g_i^{-1}g_j$. Elementary gradings arise from
the gradings on vector spaces. Let $V=\Sp{v_1,\ld, v_n}$ be a graded vector space and
$\{v_1,\ld, v_n\}$ is a graded basis such that $\deg v_i=g_i^{-1}$. Then any $E_{ij}$ is
a homogeneous linear transformation of $V$ and $\deg E_{ij}=g_i^{-1}g_j$. Any permutation
$v_i\mapsto v_{\sigma(i)}$ of basis elements induces a graded automorphism of $M_n=\mathrm{End}\:V$
and the corresponding permutation on the $n$-tuple $\tau=(g_1,\ldots,g_n)$. Hence we may
permute the components of $\tau$. Now suppose $\tau$ has the form
$$
\tau = (\underbrace{t_1,\ld,t_1}_{p_1},\ld, \underbrace{t_m,\ld,t_m}_{p_m})
$$
with $t_1,\ld, t_m$ pairwise  distinct. In this case the identity component $R_e$ is
isomorphic to $A_1\oplus\cd\oplus A_m$ where $A_i\cong M_{p_i}$, for any $i=1,\ldots,m$ and consists of all
block-diagonal matrices

$$
X=\mathrm{diag}\{ X_1,X_2,\ldots, X_m\}
$$
where $X_j$ is a $p_j\times p_j$-matrix. Moreover, for any $i\ne j$ the subspace $A_iRA_j$
is graded and all $X\in A_iRA_j$ are of degree $t_i^{-1}t_j$ in the $G$-grading. As an easy
consequence of this realization we obtain

\begin{lemma}\label{LL1}
Let $R=M_n=\oplus_{g\in G} R_g$ be a matrix algebra with an elementary $G$-grading. If
$R_e$ is simple then the grading is trivial. If $R_e$ is the sum of two simple
components, $R_e=A_1\oplus A_2$, then there exists $g\in G$, $g\neq e$, such that $A_1RA_2\subseteq
R_g$.\qed
\end{lemma}

Now we consider a matrix algebra $R=M_n$ with an involution $*: R\rightarrow R$
preserving $R_e$. Permuting $t_1,\ld, t_m$ in $\tau$ we may assume that for any $1\le j\le m-1$ either $A_j^*=A_j$ or
$A_j^*=A_{j+1}, A_{j+1}^*=A_j$. In the first case $A_j$ is simple and $A_jRA_j=A_j$. In the second case $B=A_j\oplus A_{j+1}$ is not simple
but $\ast$-simple and $A_j\simeq A_{j+1}$, i.e. $A_j$ and $A_{j+1}$ are matrix algebras
of the same size $s$. It is convenient to consider the subalgebra $BRB$ as a subset of
all matrices
$$
\mathrm{diag}\{ 0,\ldots,0,X,0,\ldots, 0\}
$$
where $X$ is $2s\times 2s$-matrix on the respective position.

Next we consider a general $G$-graded matrix algebra $R=M_n$. According to Theorem \ref{tbsz}, $R=C\otimes D$ where $C\otimes I$ is a matrix algebra with elementary grading while $I\otimes D$ is an algebra with fine grading.

\begin{lemma}\label{LL2}
Let $R=C\otimes D=\bigoplus _{g\in G} R_g$ be a $G$-graded matrix algebra with an
elementary grading on $C$ and a fine grading on $D$. Let
$\vp: R\rightarrow R$ be an antiautomorphism on $R$ preserving $G$-grading. Let also
$\vp$ acts as an involution on the identity component $R_e$ i.e $\vp^2\vert_{_{R_e}}={\rm
Id}$. Then
\begin{itemize}
\item[1)] $C_e\otimes I$ is $\vp$-stable where $I$ is the unit of $D$ and hence $\vp$
induces an involution $*$ on $C_e$;

\item[2)] there are subalgebras $B_1,\ldots,B_k\subseteq C_e$ such that $C_e=B_1\oplus\cd
\oplus B_k$, $B_1\otimes I,\ld,B_k\otimes I$ are $\vp$-stable and all $B_1,\cd,B_k$ are
$*$-simple algebras;

\item[3)] $\vp$ acts on $R_e=C_e\otimes I$ as $\vp\ast X=S^{-1}\; \tr{X} S$ where
$S=S_1\otimes I+\cd + S_k\otimes I$, $S_i\in B_iCB_i$ and $S_i=I_{p_i}$ if $B_i$ is
$p_i\times p_i$-matrix algebra with transpose involution, $S_i=\begin{pmatrix} 0 & I_{p_i} \cr
-I_{p_i} & 0\end{pmatrix}$ if $B_i$ is $2p_i\times 2p_i$-matrix algebra with symplectic involution or
$S_i=\begin{pmatrix} 0 & I_{p_i} \cr I_{p_i} & 0 \end{pmatrix}$ if $B_i\simeq M_{p_i}\oplus M_{p_i}$.

\item[4)] the centralizer of $R_e=C_e\otimes I$ in $R$ can be decomposed as $Z_1
D_1\oplus\cd\oplus Z_k D_k$ where $D_1,\ld,D_k$ are $\vp$-stable graded
subalgebras of $R$ isomorphic to $D$ and $Z_i=Z_i'\otimes I$ where $Z_i'$ is the center
of $B_i$;

\item[5)] $D$ as a graded algebra is isomorphic to $M_2\otimes\cd\otimes M_2$ where any
factor $M_2$ has the fine $(-1)$-grading.

\end{itemize}
\end{lemma}

\pp From Theorem \ref{tbsz} it follows that the identity component $R_e$ equals to
$C_e\otimes I$. Since $R_e$ is $\vp$-stable and $\vp^2={\rm Id}$ on $R_e$, the
$\vp$-action induces an involution $*$ on $C_e$. Since $C_e$ is semisimple it is a direct
sum of $*$-simple algebras, 
\bee{esum}
C_e=B_1\oplus\cd\oplus B_k.
\ene 
Now 1), 2) and 3) follows from
the classification of involution simple algebras \cite{R}.

Denote by $e_1,\ld,e_k$ the units of $B_1,\ld, B_k$, respectively. Clearly, the
centralizer $Z$ of $C_e$ in $C$ is equal to $Z_1'\oplus\cd\oplus Z_k'$ where $Z_i'$ is
the center of $B_i$ and the centralizer of $R_e$ in $R$ coincides with $Z\otimes D= Z_1
D_1\oplus\cd\oplus Z_k D_k$ where $Z_i=Z_i'\otimes I$ and $D_i=e_i\otimes D$. Obviously
the map $e_i\otimes d\mapsto d\in D$ is an isomorphism of graded algebras. Hence for
proving 4) we only need to check that all $D_1,\ld, D_k$ are $\vp$-stable.

We fix $1\le i \le k$ and consider $R'= (e_i\otimes I)R(e_i\otimes I)=C'\otimes D$ where
$C'=e_iCe_i$ and $D_i=e_i\otimes D$ is a graded subalgebra of $R'$. Then $R'$ is
$\vp$-stable since $\vp\ast(e_i\otimes I)= e_i\otimes I$. Also $R_e'=C_e'\otimes I$ with
$C_e'=B_i$. If $B_i$ is simple then $C'=C_e'$ by Lemma \ref{LL1}. In this case $D_i$ is
$\vp$-stable since $\vp$ preserves $R_e'$ and $D_i$ is the centralizer of $R_e'$ in $R'$.

Now suppose $B_i=A_1\oplus A_2$ is the sum of two matrix algebras. First we will show  that
$C'\otimes I$ is a $\vp$-stable graded subalgebra of $R$. Denote by $f_1,f_2$ the units of
$A_1$ and $A_2$ respectively. Then $f_1,f_2\in R_e$ and $\vp$ permutes $f_1,f_2$.
Moreover, $f_1C'f_2\otimes I$ and $f_2C'f_1\otimes I$ are graded subspaces. Since
$\vp\ast f_1\otimes I=f_2\otimes I, \vp\ast f_2\otimes I=f_1\otimes I$ we have
$$
\vp\ast (f_1C'f_2\otimes I) \subseteq (f_1\otimes I)R(f_2\otimes I)=f_1C'f_2\otimes D.
$$
On the other hand, since by Lemma \ref{LL1} there is $g\in G$ such that $f_1C'f_2\subseteq C_g'$ for some $g\in G$ it follows that
$$
\vp\ast (f_1C'f_2\otimes I) \subseteq R_g.
$$
Suppose now that $x\in f_1C'f_2$, $y\in D$, $x$ and $y$ are homogeneous, $\deg x=g,\deg
y=h$ Then $\deg (x\otimes y)=g$ if and only if $h=e$ that is $y=\lambda I$, for some scalar $\lambda$. It follows that
$f_1C'f_2\otimes I$ is a $\vp$-stable subspace. Similarly, $\vp\ast f_2C'f_1\otimes I
\subseteq f_2C'f_1\otimes I$, hence $C'\otimes I$ is $\vp$-stable.

Now from the decomposition $R'=C'\otimes D$ it follows that $D_i=e_i\otimes D$ is a
$\vp$-stable graded subalgebra.

For proving 5) we remark that $D$ is isomorphic to, say, $D_1$ as a $G$-graded algebra
and $D_1$ is $\vp$-stable. So, it is enough to prove that $D_1$ is the tensor product of
several copies of $M_2$. We decompose $D_1$ as the tensor product
$$
D_1\simeq R_1\otimes\cd\otimes R_m
$$
where each $R_i$ is a matrix algebra $M_{n_i}$ with a fine $\ve_i$-grading. Recall that $H=\su
D_1=H_1\times\cd\times H_m$ where $H_i\simeq {\mathbb Z}_{n_i}\times {\mathbb
Z}_{n_i}=\su R_i$, $1\le i\le m$. Now since the $G$-grading on $D_1$ is $\vp$-stable and
$$
R_i=\bigoplus_{h\in H_i}(D_1)_h
$$
it follows that $\vp\ast R_i=R_i$. Since any antiautomorphism $\vp$ on a matrix algebra
acts as $\vp\ast X=\cj{\Phi}{\,{X}}$, we can apply Lemma \ref{lfin1}. Now the proof of our lemma is complete. \bx

In what follows we discuss the canonical form of the involution $\vp$ on the whole of $R$.
As mentioned, the $\vp$-action on $R$ is defined by
$$
\vp\ast A=\cj{\Phi}{\,\tr{A}}
$$
for some matrix $\Phi$. First let $A\in R_e$. Consider the decomposition
$C_e=B_1\oplus\cd\oplus B_k$ found in Lemma \ref{LL2}. Then $A=A_1\otimes I+\cd+
A_k\otimes I$ with $A_i\in B_i, 1\le i \le k$. By Lemma \ref{LL2} $\vp$ acts on $A$ as
$$
\vp\ast A=\cj{S}{\,\tr{A}}.
$$
Hence the matrix $\Phi S^{-1}$ commutes with $\tr{A}$ for any $A\in R_e$, that is $\Phi
S^{-1}$ is an element of the centralizer of $R_e$ in $R$. Applying Claim $4)$ of Lemma
\ref{LL2} we obtain
\begin{equation}\label{einvmix}
\Phi= S_1Y_1\otimes Q_1+\cd+ S_kY_k\otimes Q_k
\end{equation}
where $Q_i\in D,Y_i\in Z_i', 1\le i\le k$. Compute now the action of $\vp^2$ on an arbitrary
$A\in R$:
$$
\vp^2\ast A=\vp\ast(\cj{\Phi}{\,\tr{A}})=\cj{\Phi}{\,\tr{(\cj{\Phi}{\,\tr{A}})}}=
(^t\Phi^{-1}\Phi)^{-1} A (^t\Phi^{-1}\Phi)
$$

Set $P=\,\tr{\,\Phi}^{-1}\Phi$.  Note that for any $T_i,T_i'\in
B_iCB_i$ and $Q_i,Q_i'\in D$, $i=1,\ld,k$, the relation
$$
\left(\sum_i T_i\otimes Q_i\right)\left(\sum_i T_i'\otimes Q_i'\right)= \sum_i T_iT_i'\otimes Q_iQ_i'
$$
holds.

We compute the value of $P$:
\bee{+1}
P=\,\tr{\,\Phi}^{-1}\Phi=\sum_{i=1}^k\;
\tr{\,(S_iY_i)}^{-1}S_iY_i\otimes\, \tr{\,Q_i}^{-1}Q_i = \sum_i\; \tr{\,S_i}^{-1}\;\tr{\,Y_i}^{-1}S_iY_i\otimes
\, \tr{\,Q_i}^{-1}Q_i.
\ene

\begin{lemma}\label{LL4}
All $Q_i$ in (\ref{+1}) satisfy $\tr{\,Q_i}^{-1}Q_i=\pm I$.
\end{lemma}

\pp Obviously it is sufficient to prove the relation
$$
e_i\otimes \tr{\iv{Q_i}}Q_i=\pm e_i\otimes I
$$
in $D_i=e_i\otimes D$. Recall that $D_i$ is $\vp$-stable (see Lemma \ref{LL2}) and $\vp$
acts on $e_i\otimes X, X\in D$ as
$$
\vp\ast (e_i\otimes X)=\cj{\Phi}{\,\tr{({e_i\otimes X})}}=(S_iY_i)^{-1}(e_i)(S_iY_i)\otimes
\cj{Q_i}{\,\tr{X}}=e_i\otimes \cj{Q_i}{\,\tr{X}}
$$
i.e. $\vp$-action induces an aniautomorphism $e_i\otimes X \mapsto e_i\otimes
\cj{Q_i}{\,\tr{X}}$ on $D_i$. By 5) Lemma \ref{LL2} $D_i$ is the tensor product
$M_2^{(1)}\otimes\cd\otimes M_2^{(r)}$ of $2\times 2$-matrix algebras with fine grading. As in
the proof of 5) Lemma \ref{LL2} we remark that all factors are $\vp$-stable. Fix a factor
$M_2^{(j)}$ and consider the action of $\vp$ on $M_2^{(j)}$. Then
$$
\vp\ast Y=\cj{T_j}{\,\tr{\,Y}}
$$
and by Lemma \ref{lfin1} $T_j= I, X_a,X_b$ or $X_{ab}$. In particular, $\tr{\,T_j}^{-1}T_j=\pm
I_2$ where $I_2$ is $2\times 2$ identity matrix. Since $e_i\otimes \cj{Q_i}{\,\tr{\,X}} = \cj{T}{(e_i\otimes\tr{\,X})}$ for all
        $X\in D$ where $T=T_1\otimes\cd\otimes T_r$ it follows that
        $e_i\otimes Q_i=\lambda T$ for some non-zero scalar $\lambda$.
        Hence $e_i\otimes Q_i$ satisfies a similar relation $e_i\otimes \tr{Q_}i^{-1}Q_i=\pm I$. \bx

We summarize what was done in this section as follows.

\begin{proposition}\label{pordtwo} Suppose $R=M_n(F)$ is the full matrix algebra over an algebraically closed field of characteristic different from 2, graded by a finite abelian group $G$. Let $\vp$ be a $G$-graded antiautomorphism of $R$ whose restriction to the identity component $R_e$ is of order two. Then $\vp$ can be given as $\vp\ast X=\cj{\Phi}{\tr{X}}$ where
\bee{eordtwo}
\Phi=S_1Y_1\otimes Q_1+\cdots+S_kY_k\otimes Q_k
\ene
where $S_i$ and $Y_i$ are described in Lemma \ref{LL2} and each $Q_i\in e_i\otimes D$ is such that $\tr{Q_}i^{-1}Q_i=\pm I$.
\end{proposition}

\section{Involutions on group graded matrix algebras}\label{si}

In this section we preserve the notation introduced earlier except that we write $\vp\ast X=X^\ast$. Our aim is to describe involutions on group graded matrix algebras. We will start with Equation (\ref{einvmix}), in which we additionally know from Lemma \ref{LL4} that $\tr{\,Q_i}^{-1}Q_i=\pm I$. Let $g^{(p)}$ mean $\underbrace{g,\ldots.g}_{q}$. Our aim is to prove the following.

\begin{theorem}\label{tinvmix} Let $\vp: X\rightarrow \cj{\Phi}{\tr{X}}$ be an involution compatible with a grading of a matrix algebra $R$ by a finite abelian group $G$. Then, after a $G$-graded conjugation, we can reduce $\Phi$ to the form
\begin{equation}\label{einvmixfin}
\Phi= S_1\otimes X_{t_1}+\cd+ S_k\otimes X_{t_k}
\end{equation}
where $S_i$ is one of the matrices $I$, $\left(\begin{array}{cc}0&I\\I&0\end{array}\right)$, or $\left(\begin{array}{cc}0&I\\-I&0\end{array}\right)$ and each $X_{t_i}$ is a matrix spanning $D_{t_i}$, $t_i\in T$. The defining tuple of the elementary grading on $C$ should satisfy the following condition. We assume that the first $l$ of summands in (\ref{einvmixfin}) correspond to those $B_i$ in (\ref{esum}) which are simple and the remaining $k-l$ to $B_i$ which are not simple. Let the dimension of a simple $B_i$ be equal to $p_i^2$ and that of a non-simple $B_{j}$ to $2p_j^2$.  Then the defining tuple has the form 
\begin{equation}\label{etuple}
\left(g_1^{(p_1)},\ldots,g_l^{(p_l)},(g_{l+1}')^{(p_{l+1})},(g_{l+1}'')^{(p_{l+1})},\ldots,(g_k')^{(p_k)},(g_k'')^{(p_k)}\right)
\end{equation}
\begin{equation}\label{econmix}
g_1^{2}t_1=\cdots=g_l^{2}t_l=g_{l+1}'g_{l+1}''t_{l+1}=\cdots=g_k'g_k''t_k.
\end{equation}
Additionally, if $\vp$ is a transpose involution then each $S_i$ is symmetric (skew-symmetric) at the same time as $X_{t_i}$, for any $i=1,\ldots,k$. If $\vp$ is a symplectic involution, then each $S_i$ is symmetric (skew-symmetric) if and only if the respective $X_{t_i}$ is skew-symmetric (symmetric), $i=1,\ldots,k$.

Conversely, if we have a grading by a group $G$ on a matrix algebra $R$ defined by a tuple as in (\ref{etuple}), for the component $C$ with elementary grading, and by an elementary abelian 2-subgroup $T$ as the support of the component $D$ with fine grading and all of the above conditions are satisfied then (\ref{einvmixfin}) defines a graded involution on $R$.
\end{theorem}
\pp
Choose $X_u\in D_u$, $u\in T$, and consider $(I\otimes X_u)^\ast$. Since $\Sp{e_1,
\ldots,e_k}$ $\otimes D$ is invariant with respect to $\vp$ we must have $(I\otimes X_u)^\ast=\alpha_1e_1\otimes X_u+\cdots+\alpha_ke_k\otimes X_u$, for some scalars $\alpha_1,\ldots,\alpha_k$. Now by Proposition \ref{pordtwo} we have
\begin{eqnarray*}
&&(I\otimes X_u)^\ast = \cj{\Phi}{\tr{(I\otimes X_u)}}\\&=&(\iv{Y_1}\iv{S_1}\otimes \iv{Q_1}+\cdots +\iv{Y_k}\iv{S_k}\otimes \iv{Q_k})\\
&\times&(e_1\otimes\tr{X_u}+\cdots+e_k\otimes\tr{X_u})
(S_1Y_1\otimes Q_1+\cd+ S_kY_k\otimes Q_k)\\
&=&e_1\otimes \cj{Q_1}{\tr{X_u}}+\cdots+e_k\otimes \cj{Q_k}{\tr{X_u}}\\
&=&e_1\otimes\alpha_1 X_u+\cdots+e_k\otimes\alpha_k X_u.
\end{eqnarray*}
It follows then that for each $i=1,\ldots,k$ we must have  $\cj{Q_i}{\tr{X_u}}$$=\alpha_i X_u$ for all $X_u\in e_i\otimes D_u$. As a result, the mapping $X\rightarrow \cj{Q_i}{\tr{X}}$ is a graded involution of an algebra with fine grading $e_i\otimes D$. By Lemma \ref{lfin1} this mapping must have the form $X_u\mapsto \cj{X_{t_i}}{X_u}$, for some $t_i\in T$. This allows us to conclude that our matrix $\Phi$ can be chosen in the form
\begin{equation}\label{einvmixmed}
\Phi= S_1Y_1\otimes X_{t_1}+\cd+ S_kY_k\otimes X_{t_k}
\end{equation}
where each $Y_i$ is in the center of $B_i$. We have that $Y_i=\lambda_ie_i$ in every case where $B_i$ is simple and $Y_i=\xi_i e_i^\prime+\zeta_ie_i^{\prime\prime}$ in every case where $B_i$ is not simple. Here $e_i^\prime$, $e_i^{\prime\prime}$ are the identities of simple components of $B_i$ and $e_i=e_i^{\prime}+e_i^{\prime\prime}$. Also any $X_{t_i}$ is either symmetric or skew-symmetric.

Now let us check the conditions (\ref{econmix}). This is done on case-by-case basis. If $U= e_iUe_j\in C$, $1\leq i,j\leq l$ then $\deg U=\iv{g_i}g_j$. Also, using (\ref{einvmix}) we obtain that $U^\ast=e_j\iv{Y_j}\iv{S_j}\tr{U}S_iY_i e_i\otimes \iv{X_{t_j}}X_{t_i}$ which is of degree $\iv{g_j}g_it_jt_i$ (we recall that by Lemma \ref{lfin1} all elements in $T$ are of order 1 or 2). Therefore, we have an equality $g_i^2t_i=g_j^2t_j$. If $U = e_iUe_j^\prime$, $1\leq i\leq l$, $l+1\leq j\leq k$, then $\deg U=\iv{g_i}g_j^\prime$. We also have that $U^\ast= e_j^{\prime\prime}\iv{Y_j}\iv{S_j}\tr{U}S_iY_ie_i\otimes \iv{X_{t_j}}X_{t_i}$, which is of degree $\iv{(g_j^{\prime\prime})}g_it_it_j$. It follows then that $g_j^{\prime}g_j^{\prime\prime}t_j=g_i^2t_i$, also in accordance with (\ref{econmix}). Finally, if $U=e_i^{\prime}U e_j^{\prime\prime}$, $l+1\leq i,j\leq k$, then $\deg U=\iv{(g_i^{\prime})}g_j^{\prime\prime}$ while $U^\ast=e_j^{\prime}\iv{Y_j}\iv{S_j}\tr{U}S_iY_ie_i^{\prime\prime}\otimes \iv{X_{t_j}}X_{t_i}$. Therefore, $\deg U^\ast=\iv{(g_j^{\prime})}g_i^{\prime\prime}t_it_j$. It then follows that
$g_j^{\prime}g_j^{\prime\prime}t_j=g_i^{\prime}g_i^{\prime\prime}t_i$, as required. The remaining three cases are in symmetry with the previous ones and produce the same results. By the way, these calculations also show that if a mapping is given by (\ref{einv1}) where $\Phi$ is as in (\ref{einvmixfin}) satisfying (\ref{econmix}) that this mapping is $G$-graded. 

Now we need to eliminate $Y_1,\ldots,Y_k$ from the formula for $\Phi$. Recall the decomposition  $R_e=B_1\oplus\cdots\oplus B_k$ from Lemma \ref{LL2}.  Each summand in (\ref{einvmixmed}) correspond to one of subalgebras $B_i$. Notice that if we apply an inner automorphism to $R$ then $\Phi$ is changed as a matrix of a bilinear form. If this automorphism is a conjugation by a matrix $P$ with identity grading then it is an isomorphism of graded algebras. In this isomorphic copy of $R=M_n$ the matrix of the involution $\vp$ will take the form of $\Phi^\prime=\tr{P}\Phi P$. We build $P$ as $P=P_1\otimes I+\cdots+P_k\otimes I$ where $P_i\in B_iCB_i$, for each $i=1,\ldots,k$. If $B_i$ is simple then $Y_i=\xi_i I$. If $B_i$ is not simple then $Y_i=\zeta_i e_i^\prime+\xi_ie_i^{\prime\prime}$. Here $e_i^\prime$, $e_i^{\prime\prime}$ are the identities of simple components of $B_i$ and $e_i=e_i^{\prime}+e_i^{\prime\prime}$. In the matrix form, $Y_i=\twm{\zeta_iI_{p_i}}{0}{0}{\xi_iI_{p_i}}$. Also, $S_i=\twm{0}{I_{p_i}}{I_{p_i}}{0}$.

Notice that since $\vp$ is an involution, $\tr{\iv{\Phi}}\Phi=\omega I$ where $\omega=\pm 1$. In other words,$\Phi=\omega\tr{\Phi}$. Now
\begin{eqnarray*} 
\tr{\Phi}&=&Y_1\tr{S_1}\otimes\tr{X_{t_1}}+\cdots+Y_k\tr{S_k}\otimes\tr{X_{t_k}}\\
&=&Y_1\tr{S_1}\otimes\alpha(t_1,t_1)X_{t_1}+\cdots+Y_k\tr{S_k}\otimes\alpha(t_k,t_k)X_{t_k}.
\end{eqnarray*}

Let us set $P_i=\frac{1}{\sqrt{\xi_i}}e_i$. If $B_i$ is simple then 
$$
S_i^\prime=\tr{P_i}S_iY_iP_i=P_iS_iY_iP_i=S_i.
$$ 
If $B_i$ is not simple then it follows  from $\tr{\Phi}=\omega\Phi$ that $\zeta_i=\xi_i\omega\alpha(t_i,t_i)$. Then 
$$S_i^\prime=\tr{P_i}S_iY_iP_i=P_iS_iY_iP_i=\frac{1}{\xi_i}\twm{0}{\xi_iI_{p_i}}{\xi_i\omega\alpha(t_i,t_i)I_{p_i}}{0}=\twm{0}{I_{p_i}}{\omega\alpha(t_i,t_i)I_{p_i}}{0}.$$

For example, if $\vp$ is a transpose involution, that is, $\Phi$ is symmetric, then $\omega=1$ and the conjugation by $P$ as above reduces $\Phi$ to the form 
$$
\Phi=S_1^\prime\otimes X_{t_1}+\cdots+ S_k^\prime\otimes X_{t_k}
$$
with
$$
\tr{\Phi}=\tr{S_1^\prime}\otimes \alpha(t_1,t_1)X_{t_1}+\cdots+ \tr{S_k^\prime}\otimes \alpha(t_k,t_k)X_{t_k}
$$
so that, according to Remark \ref{r102}, each $S_i^\prime$ is symmetric if and only if $X_{t_i}$ is symmetric, as claimed.
If $\vp$ is a symplectic involution then $\omega=-1$ and using the same equations implies that $S_i^\prime$ is symmetric if and only if $X_{t_i}$ is skew-symmetric.

The converse in the above theorem is immediate.
\bx

Now, for the determination of the gradings on simple matrix Jordan and Lie algebras, it is important to be able to compute the sets of symmetric and skew-symmetric elements of $R=M_n(F)$ under the involution just computed. A very simple remark is as follows:
$$H(R,\ast)=\mathrm{Span}\left\{\left. A+A^\ast\right| A\mbox{ from a spanning set of }R\right\},$$ $$K(R,\ast)=\mathrm{Span}\left\{\left. A-A^\ast\right| A\mbox{ from a spanning set of }R\right\}.$$ 

If $A=e_iUe_j\otimes X_u$ then $A^\ast=e_j\iv{S_j}\tr{U}S_ie_i\otimes \tr{X_{t_j}}\tr{X_u}X_{t_i}$. If we  perform obvious calculations we obtain the sets of symmetric and skew-symmetric elements of $\vp$ in the following form

\bee{esymelts}
H(R,\ast)=\Sp{e_iUe_j\otimes X_u+e_jS_j\tr{U}S_ie_i\otimes X_{t_j}\tr{X_u}X_{t_i}}
\ene
where $1\leq i,j\leq k$, $u\in T$, and $U=e_iUe_i\in C$.

Quite similarly,
\bee{esksymelts}
K(R,\ast)=\Sp{e_iUe_j\otimes X_u-e_jS_j\tr{U}S_ie_i\otimes X_{t_j}\tr{X_u}X_{t_i}}
\ene
where $1\leq i,j\leq k$, $u\in T$, and $U=e_iUe_i\in C$. Here we simultaneously replaced $\iv{S_j}$ and $\iv{X_{t_j}}$ by $S_j$ and $X_{t_j}$ 

Incidentally, this gives a canonical form for the simple graded Jordan algebras of the types $H(M_n,\ast)$ where $\ast$ is either transpose or symplectic involution (formula (\ref{esymelts})), or a simple Lie algebra of the type $B_l$, $l\geq 2$, $C_l$, $l\geq 3$, or $D_l$, $l\geq 5$(formula (\ref{esksymelts})), of which the forms suggested in \cite{BShZ} are a particular case.

\end{document}